\documentclass{amsart}

\DeclareMathOperator{\ord}{{\rm ord}}
\DeclareMathOperator{\Frob}{{\rm Frob}}
\newcommand{\abold}{{\mathbf a}}
\newcommand{\bbold}{{\mathbf b}}
\newcommand{\cbold}{{\mathbf c}}
\newcommand{\ZZ}{\mathbb Z}
\newcommand{\QQ}{{\mathbb Q}}
\newcommand{\FF}{{\mathbb F}}
\newcommand{\rr}{{\mathbf r}}
\newcommand{\sss}{{\mathbf s}}
\newcommand{\FFbar}{{\overline{\FF}}}

\newcommand{\iso}{{\stackrel{\sim}{\rightarrow}}}
\newcommand{\AAA}{{\mathcal A}}
\newcommand{\OO}{{\mathcal O}}

\newtheorem{Lemma}[subsection]{Lemma}
\newtheorem{Proposition}[subsection]{Proposition}
\newtheorem{PropDef}[subsection]{Proposition/Definition}
\newtheorem{Corollary}[subsection]{Corollary }
\theoremstyle{remark}

\title[Note on cyclotomic Euler Systems]{A note on cyclotomic Euler 
systems and the double complex method}
\author{Greg W.~Anderson}
\address{
School of Mathematics, University of Minnesota, Minneapolis,
MN 55455, USA}
\email{gwanders@math.umn.edu}

\author{Yi Ouyang}
\address{
Department of Mathematics, University of Toronto, 100 St. George St.,
Toronto, ON M5S 3G3, Canada}
\email{youyang@math.toronto.edu}
\date{October 28, 2001; revised April 24, 2002}

\begin{document}

\begin{abstract} Let $\FF$ be a finite real abelian
extension of $\QQ$. Let $M$ be an odd positive integer. For
every squarefree positive integer $r$ the prime factors of which
are congruent to $1$ modulo $M$ and split completely  in $\FF$, the 
corresponding Kolyvagin class
 $\kappa_r\in\FF^{\times}/\FF^{\times\, M}$  satisfies a
remarkable and crucial recursion which for each prime number
$\ell$ dividing
$r$ determines the order of vanishing of
$\kappa_r$ at each place of $\FF$ above $\ell$ in terms of
$\kappa_{r/\ell}$. In this note
we give the recursion a new and universal
interpretation with the help of the double complex method
introduced by Anderson and further developed by
Das and Ouyang.
  Namely, we show that the recursion satisfied by Kolyvagin classes
is the specialization of a universal 
recursion independent of $\FF$ satisfied by universal
Kolyvagin classes in the group cohomology of the
universal ordinary distribution {\em \`{a} la} Kubert tensored with
$\ZZ/M\ZZ$. Further, we show by a method involving a variant of the
diagonal shift operation introduced by Das that certain group cohomology
classes belonging (up to sign) to a basis previously constructed by
Ouyang also satisfy the universal recursion. 
\end{abstract}

\maketitle

\section{Introduction} Let $\FF$ be a finite real abelian
extension of $\QQ$. Let $M$ be an odd positive integer. For
every squarefree positive integer $r$ the prime factors of which
are congruent to $1$ modulo $M$ and split completely  in $\FF$, the  
corresponding
Kolyvagin class
\linebreak $\kappa_r\in\FF^{\times}/\FF^{\times\, M}$  satisfies a
remarkable and crucial recursion which for each prime number
$\ell$ dividing
$r$ determines the order of vanishing of
$\kappa_r$ at each place of $\FF$ above $\ell$ in terms of
$\kappa_{r/\ell}$. See Proposition 2.4 of the Rubin
appendix to Lang's text \cite{Lang1} for
a formulation of this recursion commonly employed in
the literature. In this note we actually work with a
formulation of the recursion slightly different from
but equivalent to Rubin's formulation (see
Proposition~\ref{Proposition:Crucial1} below).

The purpose of this note is to give the recursion satisfied by 
Kolyvagin classes a new and universal interpretation with the help of
the double complex method introduced by Anderson \cite{Anderson2}
and further developed by Das \cite{Das} and Ouyang
\cite{Ouyang2}. We show that the recursion is the
specialization of a universal recursion independent
of
$\FF$ satisfied by universal Kolyvagin classes in the group cohomology of
the universal ordinary distribution {\em \`{a} la} Kubert tensored with
$\ZZ/M\ZZ$ (see Proposition~\ref{Proposition:Crucial3}). 
Further, we show by a method involving a variant of the
diagonal shift operation introduced by Das \cite{Das}
that certain group cohomology classes 
belonging (up to sign) to a basis previously constructed by Ouyang
\cite{Ouyang2} also
satisfy the universal recursion (see
Corollary~\ref{Corollary:Crucial}). Taken together, our results
show that it is possible to construct classes in
$\FF^\times/\FF^{\times M}$  satisfying a useful recursion of
Kolyvagin type by methods somewhat more conceptual than have heretofore
been employed.

It is natural to expect that results similar to those presented in this
note hold for general universal Euler systems. 
Indeed, the
many beautiful results proved in Chapter
$4$ of Rubin's  book~\cite{Rubin2} strongly suggest the existence of
a general theory of universal Kolyvagin recursions.
But since there are significant technical
difficulties to deal with before the double complex
method can be brought to bear on the general theory
of Euler systems, we can for now but affirm our hope
to generalize the results of this note. In any case,
we  hope that our point of view might prove helpful
to others in the search for new applications and
new examples of Euler systems.

\section{A brief review of cyclotomic Euler systems}

\subsection{Notation and setting}
  Let $\FF$ be a finite
real abelian extension of the field $\QQ$ of rational numbers
and let $\FFbar$ be an algebraic closure of $\FF$. Let $M$
be an odd positive integer, and for every
abelian group $A$ let the $M$-torsion
subgroup be denoted $A_M$. Let
$\rr$ be the formal product of all odd primes
$\ell\equiv 1\bmod{M}$ that split completely  in
$\FF$. In the sequel we refer to formal products
of prime numbers as {\em supernatural
numbers}. Let
$G$  be the Galois group over $\FF$ of the field generated over
$\FF$ by all roots of unity in $\FFbar$ of order dividing $\rr$.
For each prime number
$\ell$ dividing $\rr$:
\begin{itemize}
\item Let $G_\ell\subset G$ be the inertia subgroup
at some place (hence all places) above $\ell$.
\begin{itemize}
\item Note that $G_\ell$ is cyclic of order $\ell-1$.
\end{itemize}
\item Let $\sigma_\ell$ be a generator of $G_\ell$.
\item Let
$N_\ell:=\sum_{i=0}^{\ell-2}\sigma_\ell^i\in \ZZ[G_\ell]$ and
$N'_\ell:=\sum_{i=1}^{\ell-2}i\sigma_\ell^{i}\in
\ZZ[G_\ell]$.
\begin{itemize}
\item Note the crucial
identity $N'_\ell(\sigma_\ell-1)=\ell-1-N_\ell$.
\end{itemize}
\item
Let $\Frob_\ell\in G/G_\ell$ be the arithmetic Frobenius
automorphism at some place (hence all places) above
$\ell$.
\item Let $\ord_\ell$ be the normalized additive
valuation of $\QQ$ corresponding to $\ell$.
\end{itemize}
For each positive integer $r$ dividing $\rr$:
\begin{itemize}
\item Let $G_r\subset G$ be the subgroup
generated by $\bigcup_{\ell\mid r}G_\ell$.
\begin{itemize}
\item Note that the evident homomorphism $\prod_{\ell\mid
r}G_\ell\rightarrow G_r$ is bijective.
\end{itemize}
\item Let $N'_r:=\prod_{\ell\mid r} N'_\ell\in \ZZ[G_r]$.
\item Let $\FF_r$ be the extension of $\FF$ generated by the
$r^{th}$ roots of unity in $\FFbar$. Put $\FF_\rr:=\bigcup_{r\mid\rr}
\FF_r$.
\item Let $\OO_r$ be the ring of
algebraic integers of $\FF_r$. Put $\OO:=\OO_1$ and
$\OO_\rr:=\bigcup_{r\mid\rr}\OO_r$.
\end{itemize}
For each positive integer $r$ dividing $\rr$ and prime number $\ell$:
\begin{itemize}
\item Let $\OO_{r,(\ell)}$ be the localization of $\OO_r$
by the multiplicative system of elements prime to
$\ell$. Put
$\OO_{(\ell)}:=\OO_{1,(\ell)}$ and
$\OO_{\rr,(\ell)}:=
\bigcup_{r\mid\rr}\OO_{r,(\ell)}$.
\end{itemize}
We fix a collection
$$\{\xi_r\in \OO_r^\times\}_{r\mid\rr}$$
of global units such that for all positive integers $r$ dividing
$\rr$  and prime numbers
$\ell$ dividing $r$ the following relations hold:
\begin{itemize}
\item
$\xi_r^{N_\ell }=
\xi_{r/\ell}^{\Frob_\ell-1}$.
\item
$\xi_{r}\equiv\xi_{r/\ell}$
modulo the radical of the ideal of $\OO_r$ generated by $\ell$.
\end{itemize}
Such a collection $\{\xi_r\}$ is called an  {\em Euler
system}.

\begin{Lemma}\label{Lemma:Satz90Etc}
Let $r$ be any positive integer dividing $\rr$.
The sequence 
\[1\rightarrow \FF_r^\times\xrightarrow{x\mapsto x^M}
\FF_r^\times\rightarrow \FF_r^\times/\FF_r^{\times M}\rightarrow 1\]
is
exact and so is the sequence
\[ 1\rightarrow\FF^\times\xrightarrow{x\mapsto
x^M}\FF^{\times M}\rightarrow H^0\left(G_r,\FF_r^\times/\FF_r^{\times
M}\right)
\rightarrow 1\]
of $G_r$-invariants. Via the latter sequence make now the identification
$$H^0\left(G_r,\FF_r^\times/\FF_r^{\times
M}\right)=\FF^\times/\FF^{\times
M}.$$ 
We have
$$
\left(\textup{\mbox{image of $H^0\left(G_r,\OO_r^\times/\OO_r^{\times
M}\right)$ in $H^0\left(G_r,\FF_r^\times/\FF_r^{\times
M}\right)$ }}\right)\subset \OO_{(\ell)}^\times/\OO_{(\ell)}^{\times
M}$$
for all prime numbers $\ell$ not dividing $r$.
\end{Lemma}
\proof 
The first sequence is
exact because under our hypotheses the field $\FF_r$  contains no
nontrivial
$M^{th}$ roots of unity. The second sequence is exact by Satz 90.
We turn to the proof of the last assertion.
Fix
$\xi\in \OO_r^\times$ representing a class in
$H^0(G_r,\OO_r^\times/\OO_r^{\times M})$ and write
$$\xi=\alpha \beta^M\;\;\;\left(\alpha\in \FF^\times,\;\;\;\beta\in
\FF_r^{\times M}\right).$$
It suffices to verify that $\alpha$ up to a factor in
$\FF^{\times M}$ belongs to $\OO_{(\ell)}^\times$.
Since $\OO_{(\ell)}$ is a principal ideal domain, it suffices to verify
that for each prime $P$ of
$\OO$ dividing
$\ell$ the order with which $P$ divides $\alpha$ is divisible by $M$.
But the latter is obvious because any prime of $\OO$ dividing $\ell$
cannot divide $r$ and hence is unramified in $\OO_r$.
\qed

\subsection{Kolyvagin classes}\label{subsection:KolyvaginClasses}
Fix a positive integer $r$
dividing $\rr$. 
For each prime number $\ell$ dividing $r$ one has
\[ \xi_r^{N'_r(1-\sigma_\ell)}\equiv \xi_r^{N'_{r/\ell} N_\ell }
\equiv \xi_{r/\ell}^{
N'_{r/\ell}(\Frob_\ell-1)}\equiv 1
\bmod{\OO_r^{\times M}}\]
by induction on the number of prime divisors of $r$ and
hence
$$\xi_r^{N'_r}\bmod{\OO_r^{\times M}}\in
H^0\left(G_r,\OO_r^\times/\OO_r^{\times M}\right).
$$
By Lemma~\ref{Lemma:Satz90Etc} there
exists a unique class
$$\kappa_r\in \FF^\times/\FF^{\times M}$$
such that
$$\xi_r^{N'_r}\equiv \kappa_r\bmod{\FF_r^{\times M}}$$
and moreover we have
\[ (\ell,r)=1\Rightarrow
\kappa_{r}\in\OO_{(\ell)}^\times/\OO_{(\ell)}^{\times
M} \]
for all prime numbers $\ell$.
We call $\kappa_r$ the {\em Kolyvagin
class} indexed by $r$.

\subsection{The operations $\nu_\ell$, $[\cdot]_\ell$ and $\exp_\ell$}
\label{subsection:NuEll} Let a prime number $\ell$ dividing
$\rr$ be given.
Let
\[ \nu_\ell:\OO_{(\ell)}^\times/\OO_{(\ell)}^{\times M}
\rightarrow (\OO/\ell)^\times_M  \]
be the unique homomorphism such that
\[ \nu_\ell \left(x \bmod
{\OO_{(\ell)}^{\times M}}\right)\equiv
x^{\frac{\ell-1}{M}}\bmod{\ell
\OO_{(\ell)}} \]
for all $x\in \OO_{(\ell)}^\times$.
Let
$$[\cdot]_\ell:\FF^\times/\FF^{\times M}
\rightarrow\left(\begin{array}{l}
\mbox{group of fractional}\\
\mbox{$\OO_{(\ell)}$-ideals}
\end{array}\right)\otimes (\ZZ/M \ZZ)$$
be the unique homomorphism such that
$$\left[x\bmod{\FF^{\times M}}\right]_\ell =\left(\begin{array}{l}
\mbox{fractional $\OO_{(\ell)}$-ideal}\\
\mbox{generated by $x$}
\end{array}\right)
\bmod{\left(\begin{array}{l}
\mbox{$M^{th}$ powers of }\\
\mbox{fractional $\OO_{(\ell)}$-ideals}
\end{array}\right)}$$
for all $x\in \FF^\times$. 
We claim that there exists a unique isomorphism 
$$\exp_\ell:\left(\begin{array}{l}
\mbox{group of fractional}\\
\mbox{$\OO_{(\ell)}$-ideals}
\end{array}\right)\otimes(\ZZ/M\ZZ)\;\iso\;
(\OO/\ell)^\times_M$$
such that
$$\exp_\ell\left( 
\left(\begin{array}{l}
\mbox{fractional $\OO_{(\ell)}$-ideal}\\
\mbox{generated by
$x^{N_\ell}$}
\end{array}\right)\otimes (1\bmod{M})\right)
\equiv
\left(x^{1-\sigma_\ell}\right)^{\frac{\ell-1}{M}}\bmod{\sqrt{\ell
\OO_{\ell,(\ell)}}}$$ for all $x\in \FF_\ell^\times$ where 
$\sqrt{\ell\OO_{\ell,(\ell)}}$ denotes the radical of
the ideal of $\OO_{\ell,(\ell)}$ generated by $\ell$.
Now each maximal ideal of $\OO_{(\ell)}$ is totally ramified 
in $\OO_{\ell,(\ell)}$, hence
every fractional
$\OO_{(\ell)}$-ideal is generated  by $x^{N_{\ell}}$ for some $x\in
\FF^{\times}_{\ell}$ unique up to a factor in $\OO_{\ell,(\ell)}^\times$
and hence
$\exp_{\ell}$ is well defined. Upon completing the
extension $\FF_\ell/\FF$ at any place of $\FF_\ell$ above $\ell$, one
obtains a Kummer extension with Galois group $G_\ell$
and hence $\exp_\ell$ is an isomorphism. The claim is proved.

\begin{Proposition}\label{Proposition:Crucial1}
For all positive integers $r$ dividing $\rr$ and prime numbers
$\ell$ dividing $r$ the identity
\[ \exp_\ell [\kappa_r]_\ell\equiv
\sqrt[\leftroot{-3}\uproot{5}M]{\xi_r^{N'_r(\sigma_\ell-1)}}\equiv
\sqrt[\leftroot{-3}\uproot{8}M]{\xi_{r/\ell}^{N'_{r/\ell}(\ell-\Frob_\ell)}}
\equiv
\nu_\ell\kappa_{r/\ell}\bmod{\sqrt{\ell\OO_{r,(\ell)}}}
\] holds, where the $M^{th}$ roots are chosen to be the unique such
existing in
$\OO_r^\times$ and $\sqrt{\ell\OO_{r,(\ell)}}$ denotes the radical of
the ideal of $\OO_{r,(\ell)}$ generated by $\ell$. (This is a
reformulation of Proposition~\textup{2.4} of the Rubin appendix to Lang's
text
\cite{Lang1}.)
\end{Proposition}
\begin{proof}  
 Write
$$\xi_r^{N'_r}=\alpha_r\beta_r^M\;\;\;(\alpha_r\in \FF^\times,
\;\;\;\beta_r\in \FF_r^\times)$$
and
$$
\xi_{r/\ell}^{N'_{r/\ell}}=\alpha_{r/\ell}\beta_{r/\ell}^M\;\;\;
\left(\alpha_{r/\ell}\in
\OO_{(\ell)}^\times,\;\;\beta_{r/\ell}\in
\OO_{r/\ell,(\ell)}^\times\right).$$
Choose $\gamma_r\in \FF_\ell^\times$ such that 
$$\mbox{$\gamma_r^{N_\ell}$
and $\alpha_r$ generate the same fractional $\OO_{(\ell)}$-ideal.}
$$
One then has
$$\gamma_r^{\frac{\ell-1}{M}}\beta_r\in\OO_{r,(\ell)}^\times.$$
Further, one has
$$
\xi_{r}^{N'_{r}(\sigma_\ell-1)}=\beta_r^{M(\sigma_\ell-1)}=
\xi_r^{N'_{r/\ell}(\ell-1)}\beta_{r/\ell}^{M(1-\Frob_\ell)}$$
and hence
$$\beta_r^{\sigma_\ell-1}=
\left(\xi_r^{N'_{r/\ell}}\right)^{\frac{\ell-1}{M}}
\beta_{r/\ell}^{1-\Frob_\ell}$$
because there are no nontrivial $M^{th}$ roots of unity in
$\FF_r$. 
Finally, one has
$$\left(\gamma_r^{1-\sigma_\ell}\right)^{\frac{\ell-1}{M}}\equiv
\beta_r^{\sigma_\ell-1}\equiv
\alpha_{r/\ell}^{\frac{\ell-1}{M}}\beta_{r/\ell}^{\ell-\Frob_\ell}\equiv
\alpha_{r/\ell}^{\frac{\ell-1}{M}}\bmod{\sqrt{\ell \OO_{r,(\ell)}}}$$
which by the definitions proves the result.
\end{proof}

\subsection{The Kolyvagin recursion}
 We say that a system of classes
$$\left\{\lambda_r\in \FF^\times/\FF^{\times M}\right\}_{r\mid\rr}$$
indexed by the positive integers dividing $\rr$
satisfies the {\em Kolyvagin recursion}
if for all positive integers $r$ dividing $\rr$
and prime numbers $\ell$
the following hold:
\begin{itemize}
\item  $(\ell,r)=1\Rightarrow
\lambda_{r}\in\OO_{(\ell)}^\times/\OO_{(\ell)}^{\times
M}$.
\item $\ell\mid  r\Rightarrow \exp_\ell [\lambda_r]_\ell=\nu_\ell
\lambda_{r/\ell}$.
\end{itemize}
In this language Proposition~\ref{Proposition:Crucial1}
and the discussion leading up to it can be condensed to
the  assertion that the system
of Kolyvagin classes satisfies the Kolyvagin recursion.

\section{Universal constructions}

\subsection{The free abelian group $\AAA$}
For each supernatural
number
$\sss$ put
$$\frac{1}{\sss}\ZZ:=\bigcup_{s\mid\sss}\frac{1}{s}
\ZZ,$$
the union being extended over all positive integers $s$
dividing
$\sss$.
Let $\AAA$ be the free abelian group generated by the family of symbols
of the form
$$[a]\;\;\;\left(a\in \frac{1}{\rr}\ZZ/\ZZ\right).
$$ 
We equip
$\AAA$ with an action of $G$ by the rule
$$\sigma[a]=[b]\Leftrightarrow \sigma\phi(a)=\phi(b)$$
for all $a,b\in \frac{1}{\rr}\ZZ/\ZZ$, injective homomorphisms
$\phi:\frac{1}{\rr}\ZZ/\ZZ
\rightarrow \FF_\rr^\times$, and $\sigma\in G$.
For each supernatural number $\sss$ dividing $\rr$, let $\AAA(\sss)$ be
the subgroup of
$\AAA$ generated by  symbols of the form
$$[a]\;\;\;\left(a\in \frac{1}{\sss}\ZZ/\ZZ\right).$$
Note that $\AAA(\sss)$ is stable under the action of $G$. 
Note that for each prime number $\ell$ dividing $\rr$
the group $\AAA(\rr/\ell)$ can
be viewed as a
$G/G_\ell$-module.
Note that
$$\AAA=\bigcup_{r\mid \rr}\AAA(r)$$
where $r$ ranges over the positive integers dividing $\rr$.

\subsection{The universal ordinary distribution}
Given any supernatural number
$\sss$ dividing
$\rr$, let
$U_\sss$ be the quotient of
$\AAA(\sss)$ by the subgroup generated by all elements of
the form
$$[a]-\sum_{\ell b=a}[b]\;\;\;
\left(
\mbox{$\ell$: a prime number dividing $\sss$},\;\;\;a\in
\frac{\ell}{\sss}\ZZ/\ZZ\right).$$
Notice that the action of
$G$ on
$\AAA(\sss)$ descends to
$U_\sss$.  Note that for every prime number $\ell$
dividing $\rr$ the group $U_{\rr/\ell}$ can be viewed as a 
$G/G_\ell$-module. 
The map
$$(a\mapsto (\mbox{class in $U_\sss$ represented by
$[a]$})):\frac{1}{\sss}\ZZ/\ZZ\rightarrow U_\sss$$ 
is the universal example of a
{\em one-dimensional ordinary distribution} of {\em level} $\sss$ {\em
\`{a} la} Kubert.  Put 
$$U:=U_{\rr}.$$ 
Abusing language slightly, we call
$U$  the {\em universal ordinary distribution}. 
See
Kubert
\cite{Kubert1}, Lang
\cite{Lang1} or Anderson's appendix to Ouyang's paper \cite{Ouyang2} for
background on the theory of the universal ordinary distribution.
By the classical results of Kubert \cite{Kubert1}, for any
supernatural number $\sss$ dividing $\rr$, the map 
$$U_\sss\rightarrow U$$ induced
by the inclusion $\AAA(\sss)\subseteq \AAA$ is an injective
homomorphism of free abelian groups with free cokernel
and hence the induced map 
$$H^0(G,U_\sss/MU_\sss)\rightarrow 
H^0(G,U/MU)$$
is also injective.
Thus we may and we do henceforth identify $U_\sss$
(resp.~$H^0(G,U_\sss/MU_\sss)$) with a subgroup of $U$ 
(resp.~$H^0(G,U/MU)$).
Note that we have
$$U_{\sss}=\bigcup_{s\mid \sss}U_s,\;\;\;
H^0(G,U_{\sss}/MU_{\sss})=\bigcup_{s\mid \sss}H^0(G,U_s/MU_s)$$
where the index $s$ in both unions ranges over the positive
integers dividing $\sss$.

\begin{Lemma}\label{Lemma:FrobNonZeroDivisor}
For every prime number $\ell$ dividing $\rr$ 
the equation 
$$(\ell-\Frob_\ell)x=0$$ has no nonzero solution $x\in U_{\rr/\ell}$.
\end{Lemma}
\begin{proof}
Fix a solution $x\in U_{\rr/\ell}$ of the equation in question.
Choose a positive integer $r$ dividing $\rr/\ell$ such that
$x\in U_r$.
Choose 
$\phi\in G_{r}$ inducing the same automorphism of
$\AAA(r)$ as does
$\Frob_\ell$. Let
$m$ be the order of
$\phi$ in the group $G_r$. Then one has an identity
$$(\ell^m-1)x=(\ell^{m-1}+\ell^{m-2}\phi+\cdots+
\ell\phi^{m-2}+\phi^{m-1})(\ell-\phi)x=0.$$
It follows  that $x=0$ because $U_{\rr/\ell}$
is a torsion-free abelian group.
\end{proof}

\subsection{The submodule $I_\ell$}
Let a prime number $\ell$ dividing $\rr$ be given.
We define 
$$I_\ell\subset U$$
to be the subgroup generated by all elements of $U$ represented by 
expressions of the form
$$[a]-[b]\;\;\;\left(a,b\in \frac{1}{\rr}\ZZ/\ZZ,\;\;\;a-b\in
\frac{1}{\ell}\ZZ/\ZZ\right).$$ Note that since
$I_\ell$ is $G$-stable and
$$(\sigma_\ell-1)U\subset I_\ell,$$
 the quotient $U/I_\ell$ can be viewed
as a
$G/G_\ell$-module.

\subsection{The resolution $L$}\label{subsection:Resolution}
Let $L$ be the free abelian group generated by symbols
of the form
$$[a,g]\;\;\;\left(\mbox{$g$: positive integer dividing $\rr$},\;\;a\in
\frac{g}{\rr}\ZZ/\ZZ\right).$$ 
We equip the abelian group $L$ with an action of $G$ by the rule
\[
\sigma[a,g]=[a',g']\;\Leftrightarrow\;\;
\sigma\phi(a)=\phi(a')\;\mbox{and}\;g=g'
\] for all symbols
$[a,g]$ and $[a',g']$ in the canonical basis of $L$, injective
homomorphisms
$\phi:\frac{1}{\rr}\ZZ/\ZZ
\rightarrow \FF_\rr^\times$ and $\sigma\in G$. 
We equip the group $L$
with a $G$-stable grading by declaring that
$$(\mbox{degree of $[a,g]$}):= -\sum_{\ell}
\ord_\ell g=-(\mbox{number of prime divisors of $g$}).$$
Here and in analogous summations below the indices $\ell$, $\ell'$
and $p$ are understood to range over prime numbers dividing $\rr$, subject
to further restrictions as noted.  We equip
$L$ with a
$G$-equivariant differential
$d$ of degree
$1$ by the rule
$$d[a,g]:=\sum_{\ell\mid g}
(-1)^{\sum_{\ell'<\ell}\ord_{\ell'}g}
\left([a,g/\ell]-\sum_{\ell b=a}[b,g/\ell]\right).$$
Now let $\sss$ be any supernatural number dividing $\rr$.
Let $L(\sss)$ be the subgroup of $L$ generated by symbols
of the form
\[ [a,g]\;\;\;\left(g:\mbox{positive integer dividing $\sss$},\;\;a\in
\frac{g}{\sss}\ZZ/\ZZ\right). \]
Then $L(\sss)$ is a $G$-  and $d$-stable graded subgroup of $L$.
It is known that
$$H^0(L(\sss),d)=U_\sss$$
via the isomorphism induced by the $G$-equivariant mapping 
$$\left([a,1]\mapsto [a]\right):(\mbox{degree zero
component of $L(\sss)$})\rightarrow
\AAA(\sss)$$ and that 
$$H^n(L(\sss),d)=0\;\;\;\mbox{for}\;n\neq 0.$$
In other words, $(L(\sss),d)$ is a resolution of $U_\sss$ in the
category of $G$-modules. See
Anderson's appendix to Ouyang's paper~\cite{Ouyang2} for details and
further discussion.

\begin{Proposition}
\label{Proposition:FormalReductionStructure}
The sequence
\[ 0\longrightarrow U_{\rr/\ell}
\xrightarrow{\ell-\Frob_\ell}
U_{\rr/\ell}\longrightarrow
U/I_\ell\longrightarrow 0
\]  is exact where the map $U_{\rr/\ell}\rightarrow U/I_{\ell}$ 
is that induced by the inclusion $U_{\rr/\ell}\subset U$.
\end{Proposition}
\begin{proof}
Let
\[ s_\ell:L\rightarrow L(\rr/\ell) \]
be the unique homomorphism
such that
\[ s_\ell[a,g]\equiv \begin{cases}
(-1)^{\sum_{\ell'<\ell}\ord_{\ell'}g}[a,g/\ell]\ &\mbox{if
$\ell\mid  g$}\\
0\ &\mbox{otherwise}
\end{cases}
\]
for all symbols $[a,g]$ in the canonical basis of $L$.
The homomorphism $s_\ell$ is of degree $1$ and satisfies the relation
$$s_\ell d=-ds_\ell$$
as can be verified by a straightforward calculation. Now consider the
sequence
$$\Sigma:\;0\rightarrow L(\rr/\ell)\longrightarrow
L/L'\stackrel{s_\ell}\rightarrow L(\rr/\ell)\rightarrow 0$$
where $L'$ is the subgroup of $L$
generated by all elements of the form
$$[a,g]-[b,g]\;\;\;\left(g\mid  \frac{\rr}{\ell},\;\;a,b\in
\frac{g}{\rr}\ZZ/\ZZ,\;\;\;a-b\in\frac{1}{\ell}\ZZ/\ZZ\right)$$
and the map $L(\rr/\ell)\rightarrow L/L'$ is that induced
by the inclusion $L(\rr/\ell)\subset L$. It is easy to verify that 
$\Sigma$ is short exact.
Since $L'$ is a graded $d$- and $G$-stable subgroup of $L$
and $(\sigma_\ell-1)L\subset L'$,
it follows that $\Sigma$ can be viewed as a
short exact sequence of complexes of
$G/G_\ell$-modules.  Because $H^*(L(\rr/\ell),d)$ is concentrated in
degree
$0$, the long exact sequence of
$G/G_\ell$-modules deduced from
$\Sigma$ by taking $d$-cohomology has at most four nonzero terms
and after making the evident identifications takes the form
$$\dots\rightarrow 0\rightarrow 
H^{-1}(L/L',d)\rightarrow
U_{\rr/\ell}\xrightarrow{1-\Frob_\ell^{-1}\ell}
U_{\rr/\ell}\longrightarrow
U/I_\ell\rightarrow 0\rightarrow
\dots$$
where the map $U_{\rr/\ell}\rightarrow U/I_\ell$ is that induced by the
inclusion $U_{\rr/\ell}\subset U$. By
Lemma~\ref{Lemma:FrobNonZeroDivisor} 
we have 
$$H^{-1}(L/L',d)=\ker\left(U_{\rr/\ell}\xrightarrow{-\Frob_\ell^{-1}
(\ell-\Frob_\ell)}
U_{\rr/\ell}\right)=0,$$ whence the result.
\end{proof}

\begin{PropDef}\label{PropDef:DeltaDef}
For every prime number $\ell$ dividing $\rr$ there exists a unique
homomorphism
$$D_\ell:H^0(G,U/MU)\rightarrow
H^0(G,U_{\rr/\ell}/MU_{\rr/\ell})$$
such that
$$\frac{(\sigma_\ell-1)x}{M}\equiv
\frac{(\ell-\Frob_\ell)y}{M}\bmod{I_\ell}\Leftrightarrow
 D_\ell(x\bmod{MU})=y\bmod{MU_{\rr/\ell}}$$
for all
$x\in U$ representing a class in $H^0(G,U/MU)$
and
$y\in U_{\rr/\ell}$ representing a class in
$H^0(G,U_{\rr/\ell}/MU_{\rr/\ell})$.  Moreover one has
$$ D_\ell H^0(G,U_r/MU_r)\subset H^0(G,U_{r/\ell}/MU_{r/\ell})$$
for all positive integers $r$ dividing $\rr$ and divisible by $\ell$.
\end{PropDef}
\proof  Put
$$\begin{array}{rcl}
X&:=&\left\{x\in U\mid \mbox{$x$ represents a class in
$H^0(G,U/MU)$}\right\},\\
Y&:=&\left\{y\in U_{\rr/\ell}\mid (\ell-\Frob_\ell)y\in
MU_{\rr/\ell}\right\},\\
Z&:=&\left\{(x,y)\in X\times Y\left|
\frac{(\sigma_\ell-1)x}{M}\equiv\frac{(\ell-\Frob_\ell)y}{M}\bmod{I_\ell}
\right.\right\}.
\end{array}
$$
Fix a positive integer $r$ dividing $\rr$ and divisible by
$\ell$.
To prove the proposition it is enough to prove the following three
claims:
\begin{enumerate}
\item $Z\cap(MU\times Y)=MU\times MU_{\rr/\ell}$.
\item $(\sigma-1)Z\subset MU\times MU_{\rr/\ell}$ for all $\sigma\in G$.
\item For all $x\in X\cap U_r$ there exists  $y\in Y\cap U_{r/\ell}$ such
that
$(x,y)\in Z$.
\end{enumerate}

We turn to the proof of the first claim. 
Only the containment
$\subset$ requires proof; 
the containment $\supset$ is trivial. Suppose
we are given 
$(x,y)\in Z\cap(MU\times Y)$. Then
$\frac{(\ell-\Frob_\ell)y}{M}\in I_\ell\cap
U_{\rr/\ell}$ and hence by
 Proposition~\ref{Proposition:FormalReductionStructure} there exists
$z\in U_{\rr/\ell}$ such that
$(\ell-\Frob_\ell)y=M(\ell-\Frob_\ell)z$.
By Lemma~\ref{Lemma:FrobNonZeroDivisor}  it follows that
$y=Mz$. Thus the first claim is proved. The second
claim follows immediately from the first.

 We turn finally to the proof of the third
claim. Let
$$\rho_\ell:\AAA(r)\rightarrow\AAA(r/\ell)$$
be the unique homomorphism such that
$$\rho_\ell[a+b]:=[a]$$
for all $a\in \frac{\ell}{r}\ZZ/\ZZ$ and $b\in \frac{1}{\ell}\ZZ/\ZZ$.
For each prime number $p$ dividing $r$, let 
$$\gamma_p:\AAA(r/p)\rightarrow\AAA(r)$$
be the unique homomorphism such that
$$\gamma_p[a]:=[a]-\sum_{pb=a}[b]$$
for all $a\in \frac{p}{r}\ZZ/\ZZ$. 
Note that $\rho_\ell$ commutes with $\gamma_p$ for $p\neq \ell$ and that 
the composite homomorphism $\rho_\ell\gamma_\ell$ induces the
endomorphism
$(1-\Frob_\ell^{-1}\ell)$ of $\AAA(r/\ell)$.
 Choose a
lifting
$\abold\in \AAA(r)$ of $x$. By hypothesis there exists an identity
$$(\sigma_\ell-1)\abold=M\bbold+\sum_{p\mid r}\gamma_p
\bbold_p\;\;\;(\bbold\in \AAA(r),\;\; \bbold_p\in \AAA(r/p)),$$ 
and hence also an identity
$$0=M\rho_\ell\bbold
-(\ell-\Frob_\ell)(\Frob_\ell^{-1}\bbold_\ell)+\sum_{p\mid
\frac{r}{\ell}}\gamma_p\rho_\ell\bbold_p.$$ Then the element
$y\in U_{r/\ell}$ represented by
$\Frob_\ell^{-1}\bbold_\ell$ has the desired property, namely that
$(x,y)\in Z$. Thus the third claim is proved and with it the
result.
\qed

\subsection{The universal Kolyvagin recursion}
We say that a family of classes
$$\{c_r\in H^0(G,U/MU)\}_{r\mid\rr}$$
indexed by the positive integers $r$ dividing $\rr$ satisfies the {\em
universal Kolyvagin recursion} if the following
conditions hold for all
positive integers $r$ dividing $\rr$ and prime numbers $\ell$:
\begin{itemize}
\item $c_r\in
H^0(G_r,U_r/MU_r)=H^0(G,U_r/MU_r)\subset
H^0(G,U/MU)$.
\item $\ell\mid r\Rightarrow  D_\ell c_r=c_{r/\ell}$.
\end{itemize}
The terminology is
justified by the next
result.
\begin{Proposition}\label{Proposition:Crucial2} Let 
$$\xi:U\rightarrow \OO_\rr^\times$$ be any $G$-equivariant homomorphism
such that
$$\xi I_\ell\subset 1+\sqrt{\ell\OO_{\rr,(\ell)}}$$
for all primes $\ell$ dividing $\rr$ where $\sqrt{\ell\OO_{\rr,(\ell)}}$
denotes the radical of the ideal of $\OO_{\rr,(\ell)}$ generated by
$\ell$. Let
$$\kappa:H^0(G,U/MU)\rightarrow
H^0(G,\FF_\rr^\times/\FF_\rr^{\times
M})\stackrel{\textup{\mbox{\scriptsize Satz 90}}}=\FF^\times/\FF^{\times
M}$$ be the homomorphism induced by
$\xi$. Let
$$\{c_r\in H^0(G,U/MU)\}_{r\mid \rr}$$
be any system of classes satisfying the universal Kolyvagin recursion.
Then the corresponding system of classes
$$\left\{\kappa c_r\in \FF^\times/\FF^{\times M}\right\}_{r\mid \rr}$$
satisfies the
Kolyvagin recursion.
\end{Proposition}
\proof Fix a positive integer $r$ dividing $\rr$ and a prime number
$\ell$. It suffices to prove the following
two assertions:
\begin{enumerate}
\item  $(\ell,r)=1\Rightarrow
\kappa c_{r}\in\OO_{(\ell)}^\times/\OO_{(\ell)}^{\times
M}$.
\item $\ell\mid r\Rightarrow \exp_\ell [\kappa c_r]_\ell=\nu_\ell
\kappa c_{r/\ell}$.
\end{enumerate}
We have $\xi U_r\subset \OO_r^\times$
by the $G$-equivariance of $\xi$, whence assertion 1
via Lemma~\ref{Lemma:Satz90Etc}. We turn to the proof of
assertion 2. By hypothesis $\ell$ divides $r$.
Fix
$$\tilde{c}_r\in U_r,\;\;\;\tilde{c}_{r/\ell}\in U_{r/\ell}$$
representing the classes $c_r$ and $c_{r/\ell}$, respectively. Write
$$\xi \tilde{c}_r:=\alpha_r\beta_r^M\;\;\;(\alpha_r\in \FF^\times,\;\;\;
\beta_r\in \FF_r^\times)$$
and
$$\xi \tilde{c}_{r/\ell}=\alpha_{r/\ell}\beta_{r/\ell}^M\;\;\;
\left(\alpha_{r/\ell}\in \OO_{(\ell)}^\times,\;\;\;
\beta_{r/\ell}\in \OO_{r,(\ell)}^{\times}\right).$$
One then has
$$
\xi\left(\frac{(\sigma_\ell-1)\tilde{c}_r}{M}
\right)=\beta_r^{(\sigma_\ell-1)},\;\;\;
\xi\left(\frac{(\ell-\Frob_\ell)\tilde{c}_{r/\ell}}{M}\right)=
\alpha_{r/\ell}^{\frac{(\ell-1)}{M}}\beta_{r/\ell}^{\ell-\Frob_\ell}$$
since there are no nontrivial $M^{th}$ roots of unity in $\OO_r^\times$.
Choose $\gamma_r\in \FF_\ell^\times$ such that 
$$\mbox{$\gamma_r^{N_\ell}$
and $\alpha_r$ generate the same fractional $\OO_{(\ell)}$-ideal.}
$$
One then has
$$\gamma_r^{\frac{\ell-1}{M}}\beta_r\in\OO_{r,(\ell)}^\times.$$
Finally, one has
$$\left(\gamma_r^{1-\sigma_\ell}\right)^{\frac{\ell-1}{M}}\equiv
\beta_r^{\sigma_\ell-1}\equiv
\alpha_{r/\ell}^{\frac{\ell-1}{M}}\beta_{r/\ell}^{\ell-\Frob_\ell}\equiv
\alpha_{r/\ell}^{\frac{\ell-1}{M}}\bmod{\sqrt{\ell \OO_{r,(\ell)}}}$$
where the crucial middle congruence holds by
Proposition/Definition~\ref{PropDef:DeltaDef} and hypothesis.
Therefore assertion 2 holds and the proposition is proved.
\qed

\section{Comparisons}

\subsection{The universal Euler system}For each positive
integer $r$  dividing $\rr$, let
$$x_r\in U_r\subset U$$
be the class represented by
$$\left[\sum_{p\mid  r}\frac{1}{p}\right]\in \AAA(r)$$
where the interior sum is extended over all primes $p$ dividing $r$.
For all positive integers $r$ dividing $\rr$ and
prime numbers $\ell$ dividing $r$ the following
clearly hold:
\begin{itemize}
\item $N_\ell x_r=(\Frob_\ell-1) x_{r/\ell}$.
\item $x_r\equiv x_{r/\ell}\bmod{I_\ell}$.
\end{itemize}
We call the family 
$$\{ x_r\in U\}_{r\mid  \rr}$$ the 
{\em universal Euler system}.

\subsection{Recovery of the Euler system by
specialization} One can easily verify the existence of a
unique $G$-equivariant homomorphism
$$\xi:U\rightarrow \OO_\rr^\times$$ such that
$$\xi x_r=\xi_r$$
for all positive integers $r$ dividing $\rr$.
Thus the given Euler system $\{\xi_r\}$
is recovered by specialization via the homomorphism $\xi$ from the
universal  Euler system $\{x_r\}$ .
Note that
$$\xi I_\ell\subset 1+\sqrt{\ell\OO_{\rr,(\ell)}}$$ 
for all primes $\ell$ dividing $\rr$.

\subsection{Universal Kolyvagin classes}
\label{subsection:UniversalKolyvagin}
Fix a positive integer $r$ dividing $\rr$.
We claim that
$$  N'_r  x_r\in U_r\subset U$$
represents a class
$$c_r\in H^0(G_r,U_r/MU_r)=H^0(G,U_r/MU_r)\subset H^0(G,U/MU).$$
For each prime
$\ell$ dividing
$r$ one has
$$(\sigma_\ell-1)N'_r x_r
\equiv
-N_\ell  N'_{r/\ell}  x_r
\equiv -(\Frob_\ell-1)  N'_{r/\ell} x_{r/\ell}\equiv 0
\bmod{MU_r}
$$
by induction on the number of prime
divisors of $r$.  Therefore $c_r$ is indeed \linebreak 
$G_r$-invariant.
 We call
$c_r$ the {\em universal Kolyvagin class} indexed by
$r$.

\subsection{Recovery of the Kolyvagin classes by
specialization} Let
$$\kappa:H^0(G,U/MU)\rightarrow
H^0(G,\FF_\rr^\times/\FF_\rr^{\times
M})\stackrel{\textup{\mbox{\scriptsize Satz
90}}}{=}\FF^\times/\FF^{\times M}$$ be the homomorphism induced by
$\xi$.
For all positive integers $r$ dividing $\rr$ one has
$$
\xi N'_rx_r=\xi_r^{N'_r}
$$
and hence
$$\kappa c_r=\kappa_r.$$
Thus the system $\{\kappa_r\}$ of Kolyvagin classes is recovered
by specialization via the homomorphism $\kappa$ from the system $\{c_r\}$
of universal Kolyvagin classes.

\begin{Proposition}\label{Proposition:Crucial3}
The universal Kolyvagin classes satisfy the universal Kolyvagin
recursion.
\end{Proposition}
\begin{proof} Fix a positive integer $r$ dividing $\rr$.
By definition the universal Kolyvagin class $c_r$ 
is represented by $N'_rx_r\in U_r$ and hence $c_r\in H^0(G,U_r/MU_r)$.
One has an identity 
$$
\frac{(\sigma_\ell-1)N'_rx_r}{M}=
\frac{\ell-1}{M}\cdot N'_{r/\ell}(x_r-x_{r/\ell})+
\frac{(\ell-\Frob_\ell)N'_{r/\ell}x_{r/\ell}}{M}$$
and hence $ D_\ell c_r=c_{r/\ell}$
by Proposition/Definition~\ref{PropDef:DeltaDef}.
\end{proof}

\subsection{Remark}
From Propositions~\ref{Proposition:Crucial2} and
\ref{Proposition:Crucial3}
one recovers the Proposition~\ref{Proposition:Crucial1}
(the latter being a reformulation of the well known
Proposition 2.4 of the Rubin appendix to Lang's text \cite{Lang1}) by
somewhat more conceptual means. We wonder if various well
known generalizations of Proposition~\ref{Proposition:Crucial1} (we have
uppermost in mind Theorem~4.5.4 on p.~91 of Rubin's book
\cite{Rubin2}) could be analogously recovered.

\section{The action of $ D_\ell$ on the
canonical basis for
$H^0(G,U/MU)$}

\subsection{The bigraded $\ZZ[G]$-module $K$}
Let $K$ be the free abelian group on symbols of the form
$$[a,g,h]\;\;\;
\left(\begin{array}{l}
\mbox{$g$: positive integer dividing $\rr$}\\
\mbox{$h$: positive integer dividing some power of $\rr$}\\
a\in \frac{g}{\rr}\ZZ/\ZZ
\end{array}\right).$$
We equip the group $K$ with a bigrading and associated total grading by
declaring that
$$\begin{array}{ccl}
\left(\mbox{bidegree of $[a,g,h]$}\right)&:=&\displaystyle
\left(-\sum_{\ell}
\ord_{\ell} g,\sum_{\ell}\ord_{\ell} h\right),\\\\
\left(\mbox{total degree of $[a,g,h]$}\right)&:=&
\displaystyle-\sum_{\ell} \ord_{\ell}
g+\sum_{\ell}\ord_{\ell} h.
\end{array}$$
We equip the group $K$ with the unique structure of bigraded
$\ZZ[G]$-module such that
$$\sigma [a,g,h]=[a',g',h']\Leftrightarrow (\sigma
\phi(a)=\phi(b)\;\mbox{and}\;g=g'\;\mbox{and}\;h=h')$$ for all symbols
$[a,g,h]$ and
$[a',g',h']$ in the canonical basis of $K$, injective homomorphisms
$\phi:\frac{1}{\rr}\ZZ/\ZZ\rightarrow\FF_\rr^\times$
and $\sigma\in G$. 
\subsection{The differentials $d$ and $\delta$}
For
each prime number
$\ell$ dividing
$\rr$ we define a $G$-equivariant differential
$d_\ell:K\rightarrow K$ of bidegree $(1,0)$ by the rule
$$d_\ell[a,g,h]:=\left\{\begin{array}{cl}
(-1)^{\sum_{\ell'<\ell}\ord_{\ell'}gh}
\left([a,g/\ell,h]-\sum_{\ell b=a}[b,g/\ell,h]\right)&\mbox{if $\ell\mid
g$,}\\\\ 0&\mbox{otherwise,}
\end{array}\right.
$$
and a $G$-equivariant differential $\delta_\ell:K\rightarrow K$ of
bidegree
$(0,1)$ by the rule
$$\delta_\ell[a,g,h]:=
(-1)^{\ord_\ell g+\sum_{\ell'<\ell}\ord_\ell gh}\left\{\begin{array}{cl}
(1-\sigma_\ell)[a,g,h\ell]&\mbox{if $\ord_\ell h\equiv 0\bmod{2}$,}\\
N_\ell[a,g,h\ell]&\mbox{if $\ord_\ell h\equiv 1\bmod{2}$.}
\end{array}\right.
$$
One can verify by a straightforward calculation that 
any two distinct operators in the family $\{d_\ell\}\cup\{\delta_\ell\}$
anticommute.  We equip $K$ with anti-commuting differentials $d$ and
$\delta$ of bidegree $(1,0)$ and $(0,1)$, respectively by the rules
$$d[a,g,h]:=\sum_{\ell}
d_\ell[a,g,h],\;\;\;\delta[a,g,h]:=\sum_{\ell}\delta_\ell[a,g,h].$$
Since the sums above contain but finitely many nonzero
terms, in fact
$d$ and
$\delta$ are well-defined. Thus we have defined a double complex
$(K,d,\delta)$ in the category of $G$-modules. 

\subsection{Comparison with Ouyang's definitions}
\label{subsection:Twist}
We define an involutive
$G$-equivariant bigraded automorphism
$\epsilon$ of
$K$ by the rule
$$\epsilon [a,g,h]:=
(-1)^{\sum_{\ell'<\ell}(\ord_{\ell}g)\cdot(\ord_{\ell'} h)}[a,g,h].$$
By a straightforward calculation one finds that 
$$\epsilon d_\ell \epsilon [a,g,h]=\left\{\begin{array}{cl}
(-1)^{\sum_{\ell'<\ell}\ord_{\ell'}g}
\left([a,g/\ell,h]-\sum_{\ell b=a}[b,g/\ell,h]\right)&\mbox{if $\ell\mid
g$,}\\\\ 0&\mbox{otherwise,}
\end{array}\right.
$$
and
$$\epsilon \delta_\ell \epsilon[a,g,h]=
(-1)^{\sum_{\ell'}\ord_{\ell'} g+\sum_{\ell'<\ell}\ord_{\ell'}
h}\left\{\begin{array}{cl} (1-\sigma_\ell)[a,g,h\ell]&\mbox{if
$\ord_\ell h\equiv 0\bmod{2}$,}\\ N_\ell[a,g,h\ell]&\mbox{if $\ord_\ell
h\equiv 1\bmod{2}$.}
\end{array}\right.
$$
It follows that $d$ (resp.~$\delta$) as defined in this paper
is $\epsilon$-conjugate to $d$ (resp.~$\delta$)
as defined by the rule appearing on p.~14 of Ouyang's paper
\cite{Ouyang2} and hence {\em mutatis mutandis} Ouyang's theory
applies to the double complex
$(K,d,\delta)$. 

\subsection{Identification of $H^0(K/MK,d+\delta)$ with $H^0(G,U/MU)$}
For any positive integer $r$ dividing
$\rr$ let
$K(r)$ be the subgroup of $K$ generated by all symbols of the form
$$[a,g,h]\;\;\;
\left(\begin{array}{l}
\mbox{$g$: positive integer dividing $r$}\\
\mbox{$h$: positive integer dividing some power of $r$}\\
a\in \frac{g}{r}\ZZ/\ZZ
\end{array}\right).$$ 
Then $K(r)$ is $G$-, $d$-, and $\delta$-stable. 
It is explained in detail in Ouyang's paper \cite{Ouyang2} how
to make the identification
$$H^*(K(r)/M K(r),d+\delta)=H^*(G_r,U_r/MU_r).$$
For our purposes in this note it is enough simply to know that the
$G$-equivariant homomorphism
$$\left((\mbox{class represented by $[a,1,1]$})\mapsto
(\mbox{class represented by $[a]$})\right)$$
$$:(\mbox{bidegree $(0,0)$
component of $K(r)/MK(r)$})\rightarrow U_r/MU_r$$ induces an isomorphism
$$H^0(K(r)/MK(r),d+\delta)\iso H^0(G_r,U_r/MU_r)=H^0(G,U_r/MU_r).$$
Then, passing to the limit over
$r$, we find that the $G$-equivariant homomorphism
$$\left((\mbox{class represented by $[a,1,1]$})\mapsto
(\mbox{class represented by $[a]$})\right)$$
$$:(\mbox{bidegree $(0,0)$
component of $K/MK$})\rightarrow U/MU$$ induces an isomorphism
$$H^0(K/MK,d+\delta)\iso H^0(G,U/MU).$$
The latter fact can be also be verified directly by
a straightforward spectral sequence argument, the key observation being
that the subcomplexes of $K$ with fixed ordinate are direct
sums of copies of  the complex $(L,d)$ discussed in
\S\ref{subsection:Resolution}.

\subsection{The canonical basis for $H^0(G_r,U_r/MU_r)$}
Fix a positive integer $r$ dividing $\rr$.
Let $S(r)$ be the bigraded $G$-,
$d$- and
$\delta$-stable subgroup of $K(r)$ generated by all symbols $[a,g,h]$ of
the canonical basis of $K(r)$ of the form
$$[a,g,h]\;\;\;(\mbox{if $a=0$, then $g$ does not divide $h$}).
$$ 
By Proposition~5.4 on
p.20 of Ouyang's paper
\cite{Ouyang2}, which is the main technical result of that paper, the
quotient map
$$(K(r),d+\delta)\rightarrow (K(r)/S(r),d+\delta)$$ is a
quasi-isomorphism of complexes. Ouyang's result is proved by
verifying that the induced map of spectral sequences is an
isomorphism at $E_2$. Clearly the family of
symbols of the form 
$$[0,g,gh]\;\;\;\left(\begin{array}{l}
\mbox{$g$: positive integer dividing
$r$}\\
\mbox{$h$: positive integer dividing some power of $r$}
\end{array}\right)
$$
forms a graded basis for $K(r)/S(r)$; moreover, it is easy to check that
$$dK(r)+\delta K(r)\subset S(r)+MK(r).$$
The upshot is that there
exists a unique
$\ZZ/M\ZZ$-basis
$$
\{\bar{c}_g\in H^0(G,U_r/MU_r)\}_{g\mid r}
$$
indexed by the positive integers $g$
dividing
$r$ such that 
the corresponding class $\bar{c}_g$ is
represented by a
$0$-cocycle of the complex $(K(r)/MK(r),d+\delta)$
congruent modulo
$S(r)+MK(r)$ to the symbol
$[0,g,g]$.  Up to signs determined by the automorphism $\epsilon$
defined in \S\ref{subsection:Twist},
the canonical basis constructed here coincides with the canonical basis
provided by Theorem~5.5 on p.~22 of Ouyang's paper
\cite{Ouyang2}.

\subsection{The canonical basis for $H^0(G,U/MU)$}
Let $S$ be the bigraded $G$-,
$d$- and
$\delta$-stable subgroup of $K$ generated by all symbols $[a,g,h]$ of
the canonical basis of $K$ of the form
$$[a,g,h]\;\;\;(\mbox{if $a= 0$, then $g$ does not divide $h$}).
$$ 
We have
$$K=\bigcup_{r\mid \rr}K(r),\;\;\;S=\bigcup_{r\mid \rr}S(r)$$
where in both unions $r$ ranges over the positive integers dividing
$\rr$. Passing to the limit over $r$ in the obvious way, we find that
there exists a unique
$\ZZ/M\ZZ$-basis
$$\{\bar{c}_r\in H^0(G_r,U_r/MU_r)\}_{r\mid \rr}$$
for $H^0(G,U/MU)$ indexed by the positive integers $r$ dividing $\rr$
such that the corresponding
class
$\bar{c}_r$ admits representation by a
$0$-cocycle of the complex \linebreak
$(K(r)/MK(r),d+\delta)$ congruent modulo $S(r)+MK(r)$
to the symbol $[0,r,r]$. We call the family
$\{\bar{c}_r\}$ the {\em canonical basis} for
$H^0(G,U/MU)$.

\subsection{The diagonal shift operator
$\Delta_\ell$} For each prime number $\ell$ dividing
$\rr$, we define the corresponding {\em diagonal
shift} operator
$\Delta_\ell$ on $K$ of bidegree $(1,-1)$   by the
rule
$$\Delta_\ell[a,g,h]:=
\left\{\begin{array}{cl}
[a,g/\ell,h/\ell]&\mbox{if $\ell\mid g$ and $\ell\mid h$,}\\
0&\mbox{otherwise.}\\
\end{array}\right.
$$
One has
$$\Delta_\ell d_{p}=d_{p}\Delta_\ell,\;\;\;
\Delta_{\ell}\delta_{p}=\delta_{p}\Delta_\ell$$
for all prime numbers $p$ distinct from $\ell$.
One has
$$\Delta_\ell d_\ell=d_\ell \Delta_\ell=0,\;\;\;
(\delta_\ell
\Delta_\ell-\Delta_\ell\delta_\ell)K\subset MK.$$
For every positive integer $r$ dividing $\rr$ one has
$$\Delta_\ell K(r)\subset \left\{\begin{array}{cl}
K(r/\ell)&\mbox{if $\ell$ divides $r$,}\\
\{0\}&\mbox{otherwise.}
\end{array}\right.$$
The action of
$\Delta_\ell$ therefore passes to 
$$H^0(K(r)/MK(r),d+\delta)=H^0(G_r,U(r)/MU(r))$$
and in the limit to
$$H^0(K/MK,d+\delta)=H^0(G,U/MU).$$  Our definition
of the diagonal shift operation $\Delta_\ell$ is
inspired by a very similar diagonal shift operation
defined on p.~3564 of the paper of Das
\cite{Das} and exploited there to great advantage.

\begin{Proposition}
For each prime number $\ell$ dividing $\rr$ the
endomorphism of \linebreak
$H^0(G,U/MU)$ induced by the diagonal shift
operation $\Delta_\ell$ coincides with
$D_\ell$. \end{Proposition}
\proof Fix a positive integer $r$ dividing $\rr$
and divisible by $\ell$.
Fix a class
$$c\in H^0(G_r,U_r/MU_r).$$
It suffices to show that $D_\ell$ and 
the endomorphism of $H^0(G,U/MU)$ induced by $\Delta_\ell$
applied to $c$ give the same result.
Let
$\cbold$ be a
$0$-chain in $K(r)$ reducing modulo $MK(r)$ to a
$0$-cycle representing $c$. Write
$$0=(d+\delta)\cbold+M\bbold$$
where $\bbold$ is a $1$-chain of $K(r)$.
For any positive integer $g$ dividing $r$
and positive integer $h$ dividing some power of
$r$ let 
$$(\abold\mapsto \abold\otimes [g,h]):\AAA(r/g)
\rightarrow K(r)$$
be the unique homomorphism such that
$$[a]\otimes [g,h]:=[a,g,h]$$
for all $a\in \frac{g}{r}\ZZ/\ZZ$. 
Write
$$\cbold=\sum \cbold_{g,h}\otimes
[g,h],\;\;\;\Delta_\ell\cbold=\sum
\cbold_{g\ell,h\ell}\otimes
[g,h]\;\;\;(\cbold_{g,h}\in
\AAA(r/g))
$$
and
$$\bbold=\sum\bbold_{g,h}\otimes[g,h]\;\;\;
(\bbold_{g,h}\in \AAA(r/g))$$
where all the sums are extended over pairs $(g,h)$
consisting of a positive integer $g$ dividing $r$
and a positive integer $h$ dividing a power of $r$. Let
$\rho_\ell:\AAA(r)\rightarrow\AAA(r/\ell)$
and $\gamma_p:\AAA(r/p)\rightarrow\AAA(r)$
be as in the proof of
Proposition/Definition~\ref{PropDef:DeltaDef}.
By hypothesis one has an identity
$$0=
\left(\sum_{\begin{subarray}{c}
p\mid r\\
p<\ell\end{subarray}}\gamma_p
\cbold_{p,\ell}\right)+\gamma_\ell
\cbold_{\ell,\ell}-\left(\sum_{\begin{subarray}{c}
p\mid r\\
p>\ell
\end{subarray}}\gamma_p
\cbold_{p,\ell}\right)+
(1-\sigma_\ell)\cbold_{1,1}+M\bbold_{1,\ell}$$
and hence also an identity
$$0=
\left(\sum_{\begin{subarray}{c}
p\mid r\\
p<\ell\end{subarray}}\gamma_p\rho_\ell
\cbold_{p,\ell}\right)-(\ell-\Frob_\ell)(\Frob_\ell^{-1}
\cbold_{\ell,\ell})-\left(\sum_{\begin{subarray}{c}
p\mid r\\
p>\ell
\end{subarray}}\gamma_p\rho_\ell
\cbold_{p,\ell}\right)+M\rho_\ell\bbold_{1,\ell}.$$
Let $x\in U_r$ be the element represented by
$\cbold_{1,1}$ and let $y\in U_{r/\ell}$ be the
element represented by
$\cbold_{\ell,\ell}$. 
One the one hand, the class of $H^0(G_r,U_r/MU_r)$ represented by 
the $0$-cocycle $\cbold \bmod{M}$ of the complex $(K(r)/MK(r),d+\delta)$
is $x\bmod{MU_r}$ and the
class of
$H^0(G_{r/\ell},U_{r/\ell}/MU_{r/\ell})$ represented by the
$0$-cocycle
$\Delta_\ell
\cbold\bmod{M}$ of the complex
$(K(r/\ell)/MK(r/\ell),d+\delta)$
is $y\bmod{MU_{r/\ell}}$.
But on the other hand, one has
$$\frac{(\sigma_\ell-1)x}{M}\equiv
\frac{(\ell-\Frob_\ell)y}{M}\bmod{I_\ell}$$
and hence
$$D_\ell(x\bmod MU)\equiv y\bmod{MU_{\rr/\ell}}$$
by Proposition/Definition~\ref{PropDef:DeltaDef}. 
Therefore the results of applying $D_\ell$ and the endomorphism of
$H^0(G,U/MU)$  induced by $\Delta_\ell$ to the class
$x\bmod{MU}$ indeed coincide.
\qed

\begin{Corollary}\label{Corollary:Crucial}
The canonical basis $\{\bar{c}_r\}$ satisfies the
universal Kolyvagin recursion.
\end{Corollary}
\proof Clear.
\qed
\begin{Corollary}\label{Corollary:TheoremB}
Any system of classes $\{b_r\}$ satisfying the universal Kolyvagin
recursion and the normalization
$b_1=\bar{c}_1$ is a $\ZZ/M\ZZ$-basis of
$H^0(G,U/MU)$.
\end{Corollary}
\proof Fix a positive integer $r$ dividing $\rr$
and arbitrarily and let 
$$r=\ell_1\cdots\ell_n$$ be the
prime factorization of $r$. One then has
$$D_{\ell_1}\cdots D_{\ell_n}b_r=b_1=\bar{c}_1
=D_{\ell_1}\cdots D_{\ell_n}\bar{c}_r$$
and hence
$$b_r-\bar{c}_r\in\ker\left(H^0(G_r,U_r/MU_r)
\xrightarrow{D_{\ell_1}\cdots D_{\ell_n}}
H^0(G_r,U_r/MU_r)\right)=\bigoplus_{\begin{subarray}{c}
r'\mid r\\ r'\neq r
\end{subarray}}\ZZ/M\ZZ\cdot \bar{c}_{r'},$$
whence the result.
\qed

\subsection{Remark}
From Corollary~\ref{Corollary:TheoremB}
it follows in particular that the system
$\{c_r\}$ of universal Kolyvagin classes defined
in \S\ref{subsection:UniversalKolyvagin} is a
$\ZZ/M\ZZ$-basis for
$H^0(G,U/MU)$. More precisely, it follows that
for every positive integer $r$ dividing $\rr$,
the family $\{c_g\}_{g\mid r}$ is a $\ZZ/M\ZZ$-basis
for $H^0(G_r,U_r/MU_r)$.
The latter fact is none other than Theorem B on p.~2
of Ouyang's paper \cite{Ouyang2}; our way of proving
Theorem B here is simpler than Ouyang's
original method.


\begin{thebibliography}{20}
\bibitem{Anderson2}
Anderson, G.\ W., {\em A double complex for computing the
sign-cohomology of the universal ordinary distribution,} Contemp.\ Math.\
{\bf 224}(1999), 1--27.

\bibitem{Das}
Das, P., {\em Algebraic Gamma monomials and double coverings of
cyclotomic fields.} Trans.~Amer.~Math.~Soc.~{\bf 352}(2000),
3557--3594

\bibitem{Kubert1}
Kubert, D.\ S., {\em The universal ordinary distribution,}
Bull.\ Soc.\ Math.\ France \textbf{107}(1979), 179--202.
\bibitem{Lang1}
Lang, S., {\em Cyclotomic Fields I and II, combined $2^{nd}$
edition,} Grad.\ Texts
in Math.\ \textbf{121}, Springer Verlag, 1990.
\bibitem{Ouyang2}
Ouyang, Y., {\em Group cohomology of the universal ordinary
distribution,} J.\ reine angew.\ Math.\ \textbf{537}(2001) 1--32.
\bibitem{Rubin2}
Rubin, K., {\em Euler systems,}  Annals of Math.\
Studies \textbf{147}, Princeton University Press, 2000.
\end{thebibliography}
\end{document}